\documentclass{amsart}[11pt,A4]
\usepackage[margin=2.54cm]{geometry}
\usepackage[T1]{fontenc}
\usepackage[latin9]{inputenc}
\usepackage{amsthm,amsmath,amssymb,color,graphicx,url}
\usepackage{float}
\usepackage{tikz}
\usepackage{graphicx}
\usepackage{mathtools}
\usepackage{makecell}
\usepackage{lineno}

\usetikzlibrary{calc}

\newtheorem{theorem}{Theorem}
\newtheorem{lemma}[theorem]{Lemma}
\newtheorem{corollary}[theorem]{Corollary}

\newtheorem{proposition}[theorem]{Proposition}
\theoremstyle{definition}
\newtheorem{definition}[theorem]{Definition}
\newtheorem{remark}[theorem]{Remark}

\newcommand{\mP}{{\mathcal P}}

\newcommand{\Aut}{\hbox{\rm Aut}}

\newcommand{\NN}{\mathbb{N}}
\newcommand{\ZZ}{\mathbb{Z}}

\newcommand{\Cov}{\hbox{\rm Cov}}

\newcommand{\rank}{\hbox{rank}}
\newcommand{\Pos}{\hbox{Pos}}

\newcommand{\D}{{\rm D}}

\newcommand{\E}{{\rm E}}

\newcommand{\M}{{\mathcal{M}}}
\newcommand{\F}{{\mathcal{F}}}
\newcommand{\B}{{\mathcal{B}}}

\newcommand{\Mon}{{\rm Mon}}

\newcommand{\Man}{{\rm Man}}

\makeatletter
\numberwithin{equation}{section}
\numberwithin{figure}{section}
\numberwithin{theorem}{section}

\makeatother

\numberwithin{equation}{section}
\numberwithin{figure}{section}

\title{Faithful and thin non-polytopal maniplexes}

\author{Dimitri Leemans}
\email{leemans.dimitri@ulb.be}
\author{Micael Toledo}
\email{micaelalexitoledo@gmail.com}
\address{Universit\'e Libre de Bruxelles, D\'epartement de Math\'ematique, C.P.216 - Alg\`ebre et Combinatoire, Boulevard du Triomphe, 1050 Brussels, Belgium}
\thanks{The authors gratefully acknowledge financial support from the F\'ed\'eration Wallonie-Bruxelles -- Actions de Recherche Concert\'ees (ARC Advanced grant)}

\begin{document}

\begin{abstract} 


Maniplexes are coloured graphs that generalise maps on surfaces and abstract polytopes. Each maniplex uniquely defines a partially ordered set that encodes information about its structure. When this poset is an abstract polytope, we say that the associated maniplex is polytopal. Maniplexes that have two properties, called faithfulness and thinness, are completely determined by their associated poset, which is often an abstract polytope. We show that all faithful thin maniplexes of rank three are polytopal. So far only one example, of rank four, of a thin maniplex that is not polytopal was known. We construct the first infinite family of maniplexes that are faithful and thin but are non-polytopal for all ranks greater than three.
\end{abstract}
\maketitle

\section{Introduction}
An abstract polytope is a partially ordered set that generalises the face-lattice of a convex polyhedron. In particular, an abstract $n$-polytope is a flagged poset satisfying two special conditions giving it additional structure: thinness (also called the diamond property) and strong connectivity. Every abstract $n$-polytope $\mP$ uniquely defines a $n$-valent edge-coloured graph $\Gamma_{\mP}$, called its flag-graph, that completely encodes its structure. This means that the poset $\mP$ can be reconstructed from the graph $\Gamma_{\mP}$, and vice versa. The possibility of jumping freely from the poset to the graph, and back, offers some obvious advantages, as it allows us to use graph theoretical tools to tackle problems pertaining to abstract polytopes (see for instance \cite{Pel2008} for a nice example of this). 

A maniplex of rank $n$ is an $n$-valent edge-coloured graph with certain properties that generalises the notion of the flag-graph of an abstract $n$-polytope. Maniplexes were first introduced in 2012 by Wilson \cite{manisteve} and have since gained considerable attention, even though an equivalent notion was previously used by Vince in \cite{vince} (where they are called linear combinatorial maps).

Given a maniplex $\M$, one can construct, in a natural way, a flagged poset that carries some information about the structure of $\M$ (see \cite{polymani}). Sometimes, but not always, this flagged poset is an abstract polytope, in which case we say $\M$ is a polytopal maniplex. A maniplex relates to its flagged poset in a similar way that the flag-graph $\Gamma_{\mP}$ relates to its associated polytope $\mP$. 
In fact, if $\M$ is isomorphic to the flag graph of some polytope $\mP$, then its associated poset $\Pos(\M)$ is a polytope isomorphic to $\mP$. In this case $\M$ is completely determined by $\Pos(\M)$. Unfortunately, this does not hold for maniplexes in general, as not all maniplexes can be reconstructed from their flagged poset. Often, some information about a maniplex $\M$ is lost in translation when we make the jump from $\M$ to $\Pos(\M)$. To ensure that enough information to fully determine $\M$ carries over to $\Pos(\M)$, certain conditions must be satisfied by $\M$. First, the nodes of $\M$ must be in one-to-one correspondence with the maximal chains of $\Pos(\M)$. This allows us to retrieve the set of flags (nodes) of $\M$ from $\Pos(\M)$. A maniplex with this property is called faithful. Secondly, it is necessary that $\M$ (or rather $\Pos(\M)$) is thin to guarantee that the edges of $\M$ are completely and uniquely determined by $\Pos(\M)$. Maniplexes that are both thin and faithful constitute the class of maniplexes that are uniquely determined by their associated flagged poset, and can be reconstructed from it.


A thin faithful maniplex is, in a sense, `almost' polytopal, as its associated poset possesses one of the two properties necessary for a flagged poset to be a polytope. It is then natural to ask whether or not there are maniplexes that are faithful and thin but not polytopal. As it transpires, no such maniplex exists for ranks smaller than $4$ (see Lemma~\ref{lem:norank3}) but an example of a rank $4$ maniplex with these properties was given in \cite{polymani}. In this paper we show that there are infinitely many thin, faithful non-polytopal maniplexes for every rank $n \geq 4$, even when we impose strong symmetry conditions. 

\begin{theorem}
\label{the:main}
For $n >3$, there are infinitely many regular maniplexes of rank $n$ that are faithful and thin but are non-polytopal. 
\end{theorem}

To prove Theorem \ref{the:main} we first prove the particular case where $n = 4$, in Proposition \ref{prop:rank4}. In fact, most of the paper will be devoted to proving the rank four case. The result extends easily to higher ranks. We make extensive use of the theory of graph covers to prove Proposition \ref{prop:rank4}. The desired infinite family of rank $4$ maniplexes is obtained by constructing maniplex covers of a particular polytopal regular $4$-maniplex that we will call $\B$. We do this in two steps: first, we construct an infinite family of (polytopal) $4$-maniplexes by constructing regular covers of $\B$, and then, we construct a special double cover for each member of this newly constructed infinite family.

The paper is structured as follows. In Section \ref{sec:pre} we give formal definitions of the notions discussed in the introduction and we prove some basic results on faithful maniplexes of small rank. In Section \ref{sec:volt} we define voltage graphs and covering maniplexes and provide a result (Lemma \ref{lem:canCIP}) that will be key to our construction. In Section \ref{sec:B}, we introduce the polytopal $4$-maniplex that will be the base graph for the construction in the following section. In Section \ref{sec:covers}, we construct an infinite family of regular $4$-maniplexes that are thin and faithful but non-polytopal. In Section \ref{sec:higher} we show that we can extend the results of the previous section to ranks greater than four.

\section{Preliminaries}
\label{sec:pre}

Before we start, let us establish some conventions and notation. All graphs in this paper are finite and simple. 
For a natural number $n$, we denote by $[n]$ the set $\{0,\ldots,n-1\}$. We use exponential notation for the action of graph automorphisms on vertices, edges or darts, so $x^\phi$ denotes the image of $x$ under a graph automorphism $\phi$, but we use standard function notation for other mappings.

\subsection{Flagged posets}

A poset $\mathcal P$ is said to be a {\em flagged poset} of rank $n$ if it has unique minimal and maximal elements and every maximal chain has length $n+2$. Clearly, every flagged poset admits a unique strictly monotone function $rank\colon \mathcal P \to \{-1,0,\ldots,n\}$, and every maximal chain contains exactly one element of rank $i$ for every $i \in \{-1,0,\ldots,n-1\}$.

We say $\mathcal P$ is {\em thin} if it satisfies the following condition:

\begin{itemize}
\item[(P.1)] For every $i \in [n]$ and for every two elements $F_{i-1}$ and $F_{i+1}$ of rank $i-1$ and $i+1$, respectively, there exist exactly two elements of rank $i$, $F_i$ and $F_i'$, such that $F_{i-1} \leq F_i, F_i' < F_{i+1}$.
\end{itemize}

The property (P.1) is usually called the diamond condition in the context of polytopes, but we will call it thinness (as such posets are called in the theory of incidence geometries) for the sake of simplicity.

Two maximal chains of $\mathcal P$ are said to be $i$-adjacent if they differ only in their element of rank $i$. Observe that in the case where the flagged poset $\mathcal P$ is thin, each maximal chain is $i$-adjacent to a unique maximal chain, and $i$-adjacency is a symmetric relation. We say $\mathcal P$ is {\em strongly connected} if it has the following property:

\begin{itemize}
\item[(P.2)] For any two maximal chains $\Phi$ and $\Psi$ of $\mathcal P$, there exists a sequence $\Phi = \Phi_0, \Phi_1, \ldots, \Phi_k = \Psi$ of maximal chains such that $\Phi \cap \Psi \subset \Phi_j$ for all $j$, and each $\Phi_j$ is $i_j$-adjacent to $\Phi_{j+1}$ for some $i_j \in [n]$.
\end{itemize}

\begin{definition}
An {\em abstract $n$-polytope} is a flagged poset of rank $n$ that is thin (property (P.1)) and strongly connected (property (P.2)).
\end{definition}

\subsection{Maniplexes}

\begin{definition}
A {\em maniplex of rank} $n$, or an {\em $n$-maniplex}, is an $n$-valent connected graph $\M$ with a proper edge-colouring $c\colon \E(\M) \to [n]$ such that if two colours $i$ and $j$ satisfy $|i-j|>1$, then every connected subgraph induced by edges with colours $i$ and $j$ is a $4$-cycle.
\end{definition} 

To avoid ambiguity, and to be more consistent with the terminology of abstract polytopes, the vertices of an $n$-maniplex are called {\em flags}, and its edges {\em links}. We  call a link of colour $i$ an {\em $i$-link} and we say that two flags connected by an $i$-link are {\em $i$-adjacent}, or that they are {\em $i$-neighbours}. Note that every flag has exactly one $i$-neighbour for each $i \in [n]$.


If $J \subset [n]$, we let $\M_J$ denote the subgraph of $\M$ induced by all the $i$-links with $i \in J$. We let $\overline{J}$ denote the complement of $J$ in $[n]$ when the context is clear. Similarly, we let $\bar{\imath}$ denote the complement of the singleton $\{i\}$. For a flag $u$ we let $\M_J(u)$ denote the connected component of $\M_J$ containing $u$. 

A {\em face} of rank $i \in [n]$, or an $i$-face, of $\M$ is a connected component of $\M_{\bar{\imath}}$, the subgraph of $\M$ resulting from the deletion of all the $i$-links. An $i$-face and a $j$-face are incident if their intersection is non-empty. 
As in abstract polytopes theory, we endow each $n$-maniplex with two {\em improper faces} of rank $-1$ and rank $n$, that are incident to all other faces of $\M$. We can think of each improper face of $\M$ as a formal copy of $\M$. The improper faces do not usually play an important role, as they do not provide any information on the structure of $\M$, but they are necessary for several notions pertaining to maniplexes and abstract polytopes to be well-defined. The faces of ranks $0$, $1$ and $(n-1)$ of an $n$-maniplex are often called {\em vertices}, {\em edges} and {\em facets}, respectively.

An automorphism of a maniplex is a graph automorphism that preserves the colouring. The set of all automorphisms of a maniplex, together with composition, forms a group that we denote $\Aut(\M)$. Since $\M$ is connected and every flag is incident to exactly one link of each colour, we see that the action of $\Aut(\M)$ on the flags is semi-regular. If $\Aut(\M)$ is transitive on the flags, then its action is regular, and we say that $\M$ is {\em regular}. 

The {\em dual} of $\M$ is the maniplex $\M^*$ with the same set of flags and links than $\M$ but with the reverse colouring. That is, a link in $\M^*$ has colour $n - c -1$ where $c$ is the colour of the corresponding link in $\M$ and $n$ is the rank of both $\M$ and $\M^*$. 
We say $\M$ is {\em self-dual} if it is isomorphic to $\M^*$.

For each $i \in [n]$, we can define an involutory permutation $r_i$ on the set of flags of $\M$, that maps each flag $u$ to it $i$-neighbour $u^i$. The group generated by these $n$ involutions in called the {\em monodromy group} of $\M$, and is denoted by $\Mon(\M)$. The group of automorphism of $\M$ can be viewed as the group of all permutations on the flags of $\M$ that commute with $\Mon(\M)$.

A walk in $\M$ is a sequence of flags $W = (u_0,u_1,\ldots,u_{k})$  such that $u_iu_{i+1}$ is a link for all $i \in [k]$. We say that $W$ is a $c_0c_1\ldots c_\ell$-walk if it traces links of colours $c_0, c_1, \ldots, c_\ell$, in that order. Note that once an initial vertex is specified, this sequence of colours uniquely determines $W$.

\subsection{Polytopal maniplexes}

Let $\M$ be an $n$-maniplex and let $F_\M$ denote the set of faces of $\M$. Let $\leq$ be a the relation on $F_\M$ given by 
\begin{align*}
F \leq G \text{ if and only if } F\cap G \neq \emptyset \text{ and } \rank(F) \leq \rank(G).
\end{align*}

As it transpires, $\leq$  is a partial order for $F_M$ and the poset $(F_\M,\leq)$ is a flagged poset of rank $n$ (see \cite[Proposition 3.1]{polymani}), which we denote $\Pos(\M)$. We say $\M$ is {\em thin} if $\Pos(\M)$ is thin. 



Conversely, if $\mP$ is a thin flagged poset of rank $n$, one can define an $n$-maniplex $\Man(\mP)$ by taking the maximal chains of $\mP$ as the flag-set of $\Man(\mP)$, and making two flags $i$-adjacent if the corresponding chains are $i$-adjacent in $\mP$. Note that the thinness of $\mP$ guarantees that $\Man(\mP)$ is well-defined. If $\mP$ is an abstract polytope, then $\Man(\mP)$ is its flag-graph. 

At first glance it may seem like the operators $\Man$ and $\Pos$ are the inverse of one another. However, this is not the case. While $\Pos(\Man(\mP)) \cong \mP$ for all thin flagged posets $\mP$, there is no guarantee that a thin maniplex $\M$ will be isomorphic to $\Man(\Pos(\M))$. Indeed, it could be that the number of maximal chains in $\Pos(\M)$ is strictly smaller than the number of flags in $\M$. In such a case, there is no possibility for $\Man(\Pos(\M))$ to be isomorphic to $\M$. An additional property, called faithfulness in \cite{polymani} must be satisfied by $\M$ for its flags to be in one-to-one correspondence with the maximal chains of $\Pos(\M)$.


\begin{definition}
An $n$-maniplex $\M$ is {\em faithful} if for every flag $\Phi$ the intersection of all faces containing $\Phi$ has only one element.
\end{definition}

When $\M$ is faithful the function $f$ mapping a maximal chain $\{F_{-1},F_0,\ldots,F_n\}$ of $\Pos(\M)$ to the unique element of $\bigcap_{i=0}^{n-1} F_i$ is well-defined and bijective. Moreover, if in addition $\M$ is thin, then $f$ maps $i$-adjacent maximal chains to $i$-adjacent flags (note that faithfulness and thinness are both necessary for this to be true in general). In this case, every path of $\M$ corresponds to a sequence of maximal chains in $\Pos(\M)$ such that any two consecutive elements of the sequences are $i$-adjacent for some $i$.

An $n$-maniplex $\M$ is {\em polytopal} if $\Pos(\M)$ is an abstract $n$-polytope. It was shown in \cite{polymani} that a faithful maniplex is polytopal if and only if it satisfied the {\em component intersection property}, defined below.



\begin{sloppypar}
\begin{definition}
An $n$-maniplex has the {\em component intersection property} (CIP) if for every set $\{H_0,H_1,\ldots,H_{k-1}\}$ of pairwise incident faces we have that $\bigcap_{j \in [k]} H_i$ is connected.
\end{definition}
\end{sloppypar}


\subsection{Faithful maniplexes of small rank}

Let us quickly prove a few auxiliary results regarding faithful maniplexes. Throughout the paper, we will be interested in whether or not a maniplex (or one of its faces) is bipartite, as bipartiteness is equivalent to orientability, in the sense that a map or an abstract polytope is orientable if and only if its associated maniplex is a bipartite graph.

\begin{lemma}
\label{lem:faithful}
Let $\M$ be a faithful $n$-maniplex and let $u$ and $v$ be two flags of $\M$. If $u$ and $v$ lie on the same $j$-face for some $j \in [n]$, then $u$ is not $j$-adjacent to $v$.
\end{lemma}

\begin{proof}
Suppose that $u$ and $v$ are $j$-adjacent. Then $F_i(u) = F_i(v)$ for all $i \neq j$ and since by hypothesis $u$ and $v$ lie on the same $j$-face, we have that $\bigcap_{i \in [n]} F_i(u) = \bigcap_{i \in [n]} F_i(v)$ and so $\{u,v\} \subseteq \bigcap_{i\in [n]} F_i(u)$, contradicting the faithfulness of $\M$. 
\end{proof}

\begin{lemma}\label{lem:norank3}
Let $\M$ be a faithful maniplex of rank $3$. If $\M$ is thin, then it is polytopal.
\end{lemma}

\begin{proof}
Since $\M$ is thin, we only need to show that $\Pos(\M)$ is strongly connected. Consider two maximal chains $\Phi$ and $\Psi$ in $\Pos(\M)$ and let $u$ and $v$ be their corresponding flags in $\M$. Let $J \subset [3]$ be the set of colours such that $\Phi$ and $\Psi$ have the same $j$-face if and only if $j \in J$. If $|J| = 2$ then $\Phi$ and $\Psi$ are $i$-adjacent (for the unique colour $i \notin J$) and there is nothing to prove. If $|J|=1$ then  $u$ and $v$ lie on the same $i$-face of $\M$ for $i \in J$. It follows that there is a $uv$-path $T$ of alternating colours $i-1$ and $i+1$ (where the addition is taken modulo $3$). Then $T$ induces a sequence of maximal chains between $\Phi$ and $\Psi$ containing $\Phi \cap \Psi$. Finally, if $J$ is empty then the connectivity of $\M$ implies the desired sequence of maximal chains exists in $\Pos(\M)$. 
\end{proof}
%
%
%
%
%

\begin{lemma}
\label{lem:edgefacebip}
If $\M$ is a faithful $4$-maniplex, then every $i$-face with $i \in \{1,2\}$ is bipartite.
\end{lemma}

\begin{proof}
Let $H$ be a $1$-face. That is, $H$ is a maximal connected subgraph of $\M$ with edges of colours in $\{0,2,3\}$. Let $C$ be a connected component of $H_{\{2,3\}}$, the subgraph of $H$ induced by edges of colours $2$ and $3$. Note that $C$ is a cycle and its edges are properly bicoloured and thus $C$ has even length $k$. Let $u_0,u_1,\ldots, u_{k-1}$ be the flags of $C$ (in cyclic order) labelled in such a way that every edge of the form $e_ie_{i+1}$ has colour $2$ if $i$ is even, and colour $3$ otherwise. Let $v_0$ be the $0$-neighbour of $u_0$. Since $\M$ is faithful, by Lemma~\ref{lem:faithful}, the flag $v_0$ is not in $C$.
Then $v_0$ belongs to another connected component $C'$ of $H_{\{2,3\}}$. Let $v_0,v_1,\ldots,v_{k'-1}$ be the flags of $C'$ labelled in such a way that edges of the form $e_ie_{i+1}$ are of colour $2$. Then $(u_1,u_0,v_0,v_1)$ is a path tracing links of colours $2$, $0$ and $2$ (that is, a $202$-path), an since $\M$ is a maniplex $v_1$ must be $0$-adjacent to $u_1$.  But now $(u_2,u_1,v_1,v_2)$ is a $303$-path, and $v_2$ must be $0$-adjacent to $u_2$. An inductive argument shows that $u_i$ is $0$-adjacent to $v_i$ for all $0 \leq i \leq k-1=k'-1$ and thus $H$ is in fact a prism of order $2k$. Since $k$ is even, $H$ is bipartite. A dual argument shows that every $2$-face is bipartite.
\end{proof}


\section{Covering graphs and maniplexes}
\label{sec:volt}

Let $\Gamma$ be a graph with vertex-set $V$ and edge-set $E$. A {\em dart} of $\Gamma$ is an ordered pair $(u,v)$ of adjacent vertices. For the sake of simplicity, and to avoid ambiguity, we will write $uv$ instead of $(u,v)$ to refer to a dart. We say that $u$ and $v$ are the {\em initial} and {\em final} vertices of $uv$, respectively. If $x := uv$ is a dart, then its {\em inverse}, denoted $x^{-1}$ is the dart $vu$. We let $\D(\Gamma)$ denote the set of darts of $\Gamma$. For a group $G$, a {\em voltage assignment} of $\Gamma$ is a function $\zeta: \D(\Gamma) \to G$ such that for every $x \in \D(\Gamma)$ we have $\zeta(x^{-1}) = \zeta(x)^{-1}$. The pair $(\Gamma,\zeta)$ is then called a {\em voltage graph}.

\begin{definition}
Let $(\Gamma,\zeta)$ be a voltage graph. The {\em derived cover} of $(\Gamma,\zeta)$, denoted $\Cov(\Gamma,\zeta)$ is the graph with vertex-set $V \times G$ where two vertices $(u,g)$ and $(v,h)$ are adjacent if $uv$ is a dart of $\Gamma$ and $h = \zeta(uv)g$.
\end{definition}

The canonical projection $\pi: \Cov(\Gamma,\zeta) \to \Gamma$ is a graph epimorphism that maps every vertex $(u,g)$ to $u$ and every edge $\{(u,g),(v,h)\}$ to $\{u,v\}$. For a vertex $u$ of $\Gamma$, the {\em fibre} of $u$ is the set $\pi^{-1}(u)$. The voltage group $G$ acts as a group of automorphisms of $\Cov(\Gamma,\zeta)$ by right multiplication on the second coordinate of every vertex. This action is semi-regular on the vertices of $\Gamma$, and transitive on each fibre. If $\Delta$ is a connected subgraph of $\Gamma$ then $G$ acts transitively on the set of connected components of $\pi^{-1}(\Delta)$ and thus all connected components of $\pi^{-1}(\Delta)$ must be isomorphic. 

\begin{remark}
\label{rem:iso}
In the particular case when $\pi^{-1}(\Delta)$ has $|G|$ connected components, all connected components must be isomorphic to $\Delta$ and no two vertices on the same component can belong to the same fibre. Indeed, in this case, the restriction of the projection $\pi$ to a connected component $H$ of $\pi^{-1}(\Delta)$ is a graph isomorphism between $H$ and $\Delta$.
\end{remark}

If $\M$ is a maniplex and $\zeta$ is a voltage assignment, we can give a proper edge-colouring to the covering graph $\Cov(\M,\zeta)$ by colouring each edge $e$ with the same colour $\pi(e)$ has in $\M$. We will henceforth assume, whenever we are dealing with the covering graph on a maniplex, that its edges are thus coloured. Note that this does not always yield a maniplex, as the covering graph could be disconnected, or a connected subgraph of two non-consecutive colours may not be a $4$-cycle. Hence, it will be important to determine when a covering graph is connected, and when are the paths of alternating non-consecutive colours $4$-cycles.

Connectedness can be characterised when the voltage assignment is ``nice'' enough. To be precise, let $\Gamma$ be a graph and let $T$ be a spanning tree of $\Gamma$. We say a voltage assignment $\zeta$ is {\em $T$-reduced} if it is trivial on every dart of $T$. The following lemma is folklore.


\begin{lemma}
\label{lem:voltcon}
Let $(\Gamma,\zeta)$ be a voltage graph for a voltage group $G$. If $\zeta$ is $T$-reduced for some spanning tree $T$ of $\Gamma$, then $\Cov(\Gamma,\zeta)$ is connected if and only if the images of $\zeta$ generate $G$. 
\end{lemma}

This settles the problem of connectivity. As for the second condition, let us consider a walk $W = (u_0, u_1, \ldots, u_{k})$ in $\Gamma$. If $\zeta$ is a voltage assignment for $\Gamma$, then the {\em net voltage} of $W$ is the product $\zeta(x_{k-1})\zeta(x_{k-2})\ldots\zeta(x_1)$ where $x_i$ denotes the dart $u_{i}u_{i+1}$. A classic result in the theory of graph covers states that a walk in $\Cov(\Gamma,\zeta)$ is closed if and only if $\pi(W)$ is a closed walk of $\Gamma$ with trivial net voltage. From the discussion in the above paragraphs, we obtain the following lemma.

\begin{lemma}
\label{lem:covermani}
Let $\M$ be an $n$-maniplex and let $\zeta$ be a voltage assignment for $\M$. The covering graph $\Cov(\M,\zeta)$ is an $n$-maniplex if and only if it is connected, and every closed walk in $\M$ of length four with alternating colours $i$ and $j$, $|i-j|>1$, has trivial net voltage.
\end{lemma}

When $\M$ is a $4$-maniplex with a covering graph $\Cov(\M,\zeta)$ that is also a $4$-maniplex, we can determine sufficient conditions on $\M$ for $\Cov(\M,\zeta)$ to be thin.

\begin{lemma}
\label{lem:diamond}
Let $\M$ be a $4$-maniplex and let $\zeta$ be a voltage assignment with a voltage group $G$ such that $\bar{\M} := \Cov(\M,\zeta)$ is a maniplex. Assume that for every $i$-face $F$ of $\M$, the number of connected components of $\pi^{-1}(F)$ is $1$ if $i \in \{0,3\}$ or $|G|$ if $i \in \{1,2\}$. Then, if $\M$ is thin,  so is $\Cov(\M,\zeta)$. 
\end{lemma}

\begin{proof}
Suppose $\M$ is thin. Let  $k = |G|$ and $i \in [4]$. 
Let $\bar{F}_-$ and $\bar{F}_+$ be two faces of $\bar{\M}$ of ranks $(i-1)$ and $(i+1)$ respectively and let $\bar{\Delta} = \bar{F}_- \cap \bar{F}_+$. Now, suppose for a contradiction that there are three distinct $i$-faces $H_1$, $H_2$ and $H_3$ of $\bar{\M}$ that intersect $\bar{\Delta}$. Let $F_- = \pi(\bar{F}_-)$ and $F_+ = \pi(\bar{F}_+)$, and note that $F_-$ and $F_+$ are faces of $\M$ of ranks $i-1$ and $i+1$, respectively. By hypothesis, one of $\pi^{-1}(F_-)$ and $\pi^{-1}(F_+)$ must be connected, while the other has $k$ connected components (observe that if $F$ is an improper face of $\M$, then $\pi^{-1}(F)$ is connected). Since $\pi^{-1}(F_- \cap F_+) = \pi^{-1}(F_-) \cap \pi^{-1}(F_+)$, we see that $\pi^{-1}(F_- \cap F_+)$  has $k$ connected components, one of which is $\bar{\Delta}$. Then, no two flags of $\bar{\Delta}$ belong to the same fibre. It follows that the projections of $H_1$, $H_2$ and $H_3$ are all distinct. Therefore $\pi(H_1)$, $\pi(H_2)$ and $\pi(H_3)$ are distinct $i$-faces of $\M$ that intersect $F_- \cap F_+$, contradicting that $\M$ is thin.
\end{proof}

\subsection{Canonical double covers}

A particular case of the derived cover construction, called a canonical double cover, will be very useful throughout the remainder of this paper.

\begin{definition}
The {\em canonical double cover} of a simple graph $\Gamma$ is the derived cover $\Cov(\Gamma, \zeta)$ where $\zeta: \D(\Gamma) \to \ZZ_2$ is the voltage assignment given by $\zeta(x) = 1$ for all $x \in D$, where $1$ denotes the additive non-trivial element of $\ZZ_2$.
\end{definition}

The canonical double cover of a graph has a series of important and rather useful properties, which follow easily from the definition, and that we summarize in the following remark.

\begin{remark}
\label{rem:can}
If $\bar{\Gamma}$ is the canonical double cover of a connected simple graph $\Gamma$ then:
\begin{enumerate}
\item $\bar{\Gamma}$ is a bipartite simple graph;
\item $\bar{\Gamma}$ is disconnected if and only if $\Gamma$ is bipartite;
\item if $\bar{\Gamma}$ is disconnected, it consists of two isomorphic copies of $\Gamma$;
\item every automorphism $\alpha$ of $\Gamma$ extends naturally to an automorphism $\bar{\alpha}$ of $\bar{\Gamma}$ acting by the rule $(u,g)^{\bar{\alpha}} = (u^\alpha,g+1)$;
\item $\bar{\Gamma}$ admits a group of automorphisms isomorphic to $\Aut(\Gamma) \times \ZZ_2$.
\end{enumerate}
\end{remark}



In the particular case when $\Gamma$ is a $4$-maniplex, a lot more can be said about its canonical double cover, provided that the $0$- and $3$-faces of $\Gamma$ are non-bipartite. The following lemma will be, in a sense, the key to the construction that will give us Theorem \ref{the:main}.

\begin{lemma}
\label{lem:canCIP}
Let $\M$ be a faithful $4$-maniplex such that every $i$-face with $i \in \{0,3\}$ is non-bipartite. Then, the following hold for its canonical double cover $\overline{\M}$: 
\begin{enumerate}
\item if $\M$ is regular then so is $\overline{\M}$;
\item if $\M$ is thin, then so is $\overline{\M}$;
\item if $\M$ is faithful, then so is $\overline{\M}$;
\item $\overline{\M}$ does not satisfy the CIP;
\item $\overline{\M}$ has the same number of $0$- and $3$-faces as $\M$, but twice the number of $1$- and $2$-faces.
\end{enumerate}
\end{lemma}

\begin{proof} 
Items (1) and (2) follow from Remark \ref{rem:can} and Lemma \ref{lem:diamond} respectively. We will prove item (3) last.

Let us now prove that (4) holds. For $i \in \{0,3\}$ let $F_i$ be an $i$-face of $\M$ and let $\bar{F_i} = \pi^{-1}(F_i)$. Observe that each $\bar{F_i}$ is connected and thus $\bar{F_0}$ and $\bar{F_3}$ are faces of rank $0$ and $3$ of $\overline{\M}$ respectively. Let $\Delta = F_0 \cap F_3$ and  $\bar{\Delta} = \pi^{-1}(\Delta)$, and observe that $\bar{\Delta} = \bar{F_0} \cap \bar{F_3}$. We will show that $\bar{\Delta}$ is disconnected. If $\Delta$ is disconnected, then so is $\bar{\Delta}$, and we are done. Let us assume that $\Delta$ is connected. Since $\Delta$ only has edges of colours $1$ and $2$, it must be a cycle of even length, and thus a bipartite graph. Then $\pi^{-1}(\Delta)$ is disconnected and $\overline{\M}$ does not have the CIP. 

To prove (5) note that every $i$-face of $\overline{\M}$ is a connected component of $\pi^{-1}(F)$ for some $i$-face $F$ of $\M$. Let $F$ be an $i$-face of $\M$. If $i \in \{0,3\}$, then $F$ is non-bipartite and by Remark \ref{rem:can}, $\pi^{-1}(F)$ is connected. If $i \in \{1,2\}$ then by Lemma \ref{lem:edgefacebip} the $i$-face $F$ is bipartite and thus $\pi^{-1}(F)$ has two connected components.

Finally, to prove (3) assume $\M$ is faithful. Consider a set of pair-wise incident faces $\{\bar{F_0},\bar{F_1},\bar{F_2},\bar{F_3}\}$ of $\overline{\M}$ and let $F_i := \pi(\bar{F}_i)$. Note that each $F_i$ is an $i$-face of $\M$. By item (5), $\pi^{-1}(F_1)$ has two connected components and by Remark \ref{rem:iso}, $\pi$ maps $\bigcap_{i=0}^3 \bar{F}_i$ isomorphically to $\bigcap_{i=0}^3 F_i$ which consists of a singleton, since $\M$ faithful. 
\end{proof}

\section{The $4$-polytope $\{\{4,3\}_3,\{3,4\}_3\}$}
\label{sec:B}

Our goal is to show that there are infinitely many regular $4$-maniplexes that are thin and faithful, but that do not have the CIP (and are therefore not polytopal). We know that such a maniplex can be obtained by constructing the canonical double cover of a thin, faithful regular $4$-maniplex satisfying the hypotheses of Lemma~\ref{lem:canCIP}. In this section we will briefly present a regular $4$-polytope satisfying these conditions and then, in the following section, use it to construct an infinite family of faithful thin regular $4$-maniplexes that agree with the hypotheses of Lemma~\ref{lem:canCIP}. 

\begin{figure}[h!]
\begin{center}
\includegraphics[width=0.60\textwidth]{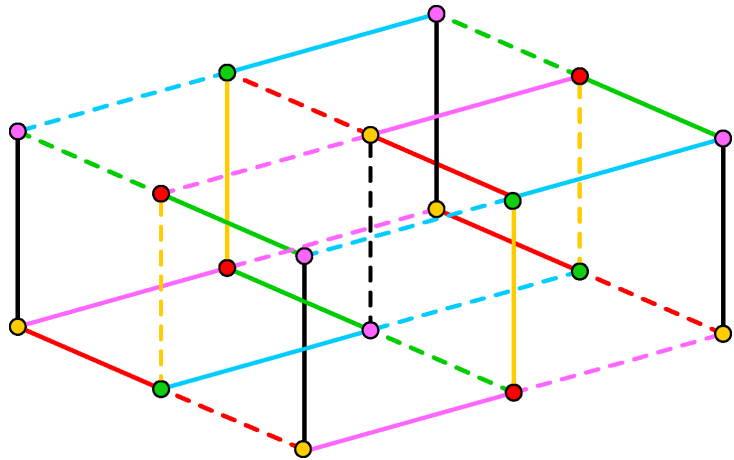}
\caption{A drawing of the skeleton of $\B$, prior to identification.}
\label{fig:skeletonB}
\end{center}
\end{figure}

The basis for our construction will be the regular $4$-polytope known as $\{\{4,3\}_3,\{3,4\}_3\}$ in the Atlas of Small Regular Abstract Polytopes \cite{atlas}. We will slightly abuse language and identify  $\{\{4,3\}_3,\{3,4\}_3\}$ with its corresponding maniplex, which we call $\B$, and give a brief description of its structure. Further details can be found in \cite{atlas}. We can construct $\B$ by glueing four copies of the hemi-cube in such a way that every two facets share exactly one (square) face, resulting in a regular $4$-polytope having $4$ vertices, $6$ edges, $6$ faces and $4$ facets. To get a better grasp of what $\B$ may look like, let us get our attention to Figure \ref{fig:skeletonB}. The drawing in the center is the "unfolded" skeleton of $\B$, prior to identification, which consists of 4 cubes arranged around a central edge (the black vertical dashed edge). We obtain the actual skeleton by identifying vertices of the same colour and edges of the same colour (regardless of whether they are solid or dashed). With this identification, each of the four cubes will collapse into one of the hemi-cubic facets of $\B$. The skeleton of $\B$ is then isomorphic to the complete graph $K_4$. The faces (of rank $2$) of $\B$ correspond to the $4$-cycles in the drawing where two cycles $C$ and $C'$ are equivalent if they trace edges of the same colours in the same order and either every edge of $C$ has the same style (solid or dashed) as its counterpart in $C'$, or every edge of $C$ has a different style as its counterpart in $C'$.

The vertex-figures of $\B$ are hemi-octahedral and $\B$ is in fact self-dual. Thus, the $0$- and $3$-faces of $\B$ are not bipartite but the $1$- and $2$-faces are (see Lemma~\ref{lem:edgefacebip}). Moreover, since $\B$ is polytopal, it is thin and faithful, and by Lemma \ref{lem:canCIP} the canonical double cover $\overline{\B}$ of $\B$ is a thin, faithful, regular $4$-maniplex but it does not have the CIP. It has $4$ facets, $12$ faces, $12$ edges and $4$ vertices. Its skeleton can also be obtained from the drawing in Figure \ref{fig:skeletonB}, by identifying vertices having same colour, and edges having both the same colour and same line style (either dashed or solid). The skeleton is a non-simple graph isomorphic to $K_4$ with every edge doubled.

An exhaustive computer search performed in {\sf Gap}~\cite{GAP4} shows that, among regular $4$-maniplexes, $\overline{\B}$ is the smallest (with fewest flags) to be non-polytopal but to be thin and faithful. Since no such maniplex exists for ranks smaller than $4$, $\overline{\B}$ is, in a sense, the smallest thin, faithful, non-polytopal, regular maniplex.

\section{The covers of $\B$}
\label{sec:covers}

As stated before, the goal of this section is to construct an infinite family of faithful thin regular $4$-polytopes that satisfy the conditions of Lemma \ref{lem:canCIP}. Then, by taking their canonical double covers, we will obtain the desired infinite family of thin, non-polytopal regular $4$-maniplexes. 

To achieve this, we define for each positive integer $n$, a voltage assignment $\zeta_n: \D(\B) \to \ZZ_n$ for $\B$ in such a way that the regular cover of $(\B,\zeta_n)$ has the desired properties. Such an assignment gives a non-trivial voltage (either $1$ or $-1$) to darts belonging to $16$ prescribed $2$-links, and voltage $0$ everywhere else.

\begin{figure}[h!]
\begin{center}
\includegraphics[width=0.80\textwidth]{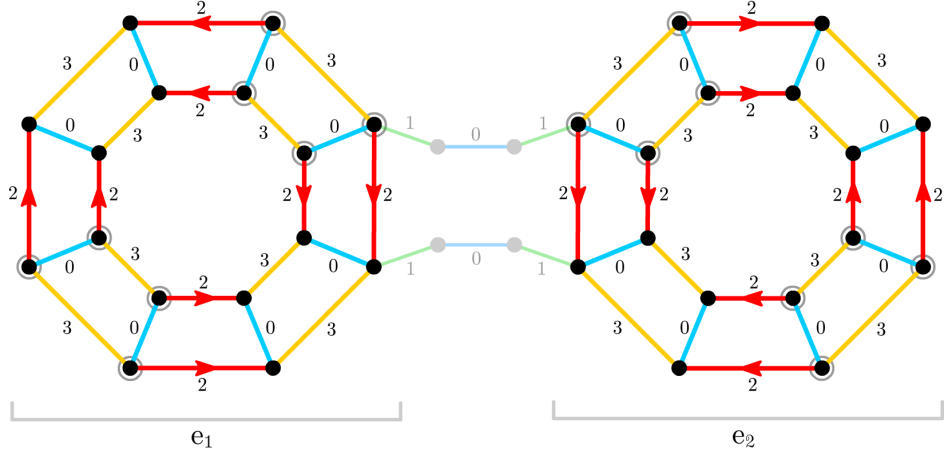}
\caption{A part of $\B$ containing the edges $e_1$ and $e_2$ where the $16$ darts belonging to $X$ are indicated by an arrowhead.}
\label{fig:edgeX}
\end{center}
\end{figure}

\begin{figure}[h!]
\begin{center}
\includegraphics[width=0.70\textwidth]{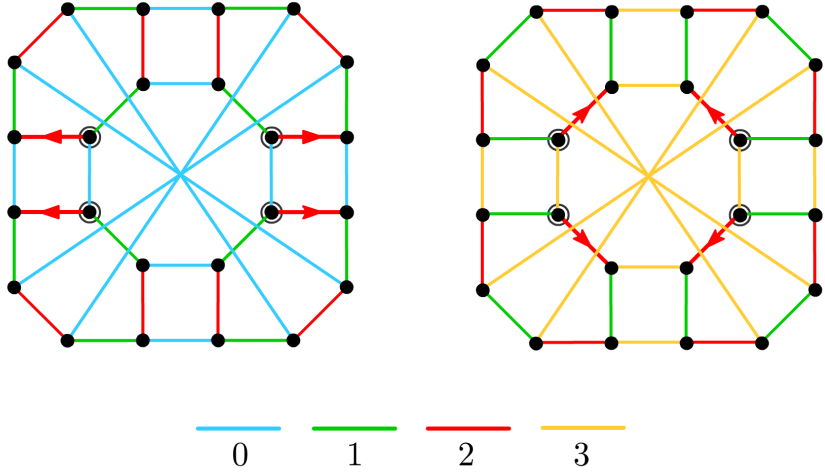}
\caption{A facet and a vertex with voltage assignment $\zeta_n$. Links without arrows have trivial voltage.}
\label{fig:vfX}
\end{center}
\end{figure}

Consider two edges ($1$-faces) $e_1$ and $e_2$ of $\B$, that do not share a vertex (for instance, the black and yellow edges in Figure \ref{fig:skeletonB}). Recall that each edge is a subgraph of $\B$ and, in particular, each of $e_1$ and $e_2$ is isomorphic to an octagonal prism with links of colours $0$, $2$ and $3$. We let $X$ be the set of darts shown in Figure \ref{fig:edgeX}, where a $2$-link with an arrowhead directed from, say, $u$ to $v$ indicates that the dart $uv$ is in $X$. Denote by $\F_X$ the set of flags incident to a dart in $X$ (flags marked with a circle in Figures \ref{fig:edgeX} and \ref{fig:vfX}). The following holds:
\begin{enumerate}
\item $u \in \F_X$ implies $u^0,u^3 \in \F_X$;
\item $u \in \F_X$ implies $u^1,u^2 \notin \F_X$;
\item $u \in \F_X$ implies $u^{101} \in \F_X$.
\end{enumerate}

Once the edges $e_1$ and $e_2$ are fixed, there are exactly two (equivalent) ways to chose $X$ so that the above properties hold: one by orienting the $2$-links like in Figure \ref{fig:edgeX}, and one by taking the opposite orientation. Each facet and each vertex of $\B$ will have exactly four darts in $X$ (see Figure \ref{fig:vfX}).

We can now define a voltage assignment for $\B$. For a positive integer $n$, we let $\zeta_n: \D(\B) \to \ZZ_n$ be given by
\[
\zeta_n(x)=
 \begin{cases}
1 & \text{ if $x \in X$} \\ 
-1 &\text{ if $x^{-1} \in X$} \\
0 & \text{ otherwise}
\end{cases}
\]

For $n \in \NN$, let $\B^n$ denote the derived cover $\Cov(\B,\zeta_n)$. Observe that $\B^n$ is a $4$-maniplex. To see this, first note that the set of links having a dart in $X$ is not a cut set of $\B$, and thus $\B^n$ is connected by Lemma \ref{lem:voltcon}. Furthermore, if $C$ is a cycle of length $4$ in $\B$ with alternating colours $i$ and $j$, $|i-j|>1$, then either all darts of $C$ have trivial voltage, or $C$ has two $2$-links whose voltages cancel each other. In either case, the net voltage of $C$ is trivial, and by Lemma \ref{lem:covermani}, $\B^n$ is a maniplex. 

\begin{lemma}
Let $n \in \NN$ and let $\zeta_n$ be the voltage assignment for $\B$ defined above. If $F_i$ is an $i$-face of $\B$, then $\pi^{-1}(F_i)$ has one connected component if $i \in \{0,3\}$ or $n$ connected components if $i \in \{1,2\}$. 
\label{lem:coverfaces}
\end{lemma}

\begin{proof}
First suppose $i \in \{0,3\}$. It is clear from the definition of $\zeta_n$ that $F_i$ admits a spanning tree $T$ whose darts all have trivial voltage (this can readily be seen in Figure \ref{fig:vfX}). That is, the restriction of $\zeta_n$ to $F_i$ is $T$-reduced and, by Lemma \ref{lem:voltcon}, $\pi^{-1}(F_i)$ is connected. 

If $i = 2$, then $F_i$ has no $2$-links, and thus no darts in $X$. By the definition of $\zeta_n$, all darts of $F_i$ have trivial voltage. It follows that $\pi^{-1}(F_i)$ consists of $n$ isomorphic copies of $F_i$.

Finally, if $i = 1$ then there are two possibilities for $F_i$. Either all darts of $F_i$ have trivial voltage (and we are done) or $F_i$ is one of the two distinguished edges $e_1$ or $e_2$. We may assume without loss of generality that $F_i$ is the edge $e_1$ and thus $F_i$, together with the corresponding restriction of $\zeta_n$, is precisely the voltage graph in the left side of Figure \ref{fig:edgeX}. Let $U$ be the set of flags of $F_i$ contained in $\F_X$, and let $V$ be the remaining eight flags (all of which are $2$-adjacent to a flag in $\F_X$). Then, for every $j \in [n]$, the subgraph $H_j$ of $\pi^{-1}(F_i)$ induced by the set of flags $\{(u,j) \mid u \in U\} \cup \{(v,j+1) \mid v \in V\}$ is a connected component of $\pi^{-1}(F_i)$. Moreover, each $H_j$ is isomorphic to $F_i$ and there are precisely $n$ distinct connected components. 
\end{proof}


\begin{corollary}
If $n$ is a positive integer, then $\B^n$ is faithful.
\end{corollary}

\begin{proof}
Use the same argument as in the proof of item (3) of Lemma \ref{lem:canCIP}.
\end{proof}

We are now ready to prove that $\B^n$ satisfies all the hypotheses of Lemma \ref{lem:canCIP}. We will prove a little bit more, and show that $\B^n$ is in fact also polytopal. This is not necessary for our construction, but it is indeed nice that the covering maniplexes $\B^n$ are polytopes whose canonical double covers are non-polytopal.

\begin{lemma}
For all $n \in \NN$, $\B^n$ is a faithful thin $4$-maniplex with the CIP and such that the $0$- and $3$-faces are non-bipartite.
\end{lemma}

\begin{proof}
Observe that, since $\B$ is polytopal, $\B^n$ is thin by Lemmas \ref{lem:coverfaces} and \ref{lem:canCIP}.
 
Let us show that $\B^n$ has the CIP. Let $\mathcal{H} = \{\bar{H}_{i_0},\ldots,\bar{H}_{i_k}\}$ be a set of pairwise incident faces of $\B^n$. We need to show that the intersection $\bigcap \mathcal{H}$ is a connected subgraph. We will divide the analysis in two cases, depending on whether or not $\mathcal{H}$ contains an $i$-face with $i \in \{1,2\}$.

First, suppose it does not. Then $\mathcal{H}$ consists of a $0$-face $\bar{H}_0$ and a $3$-face $\bar{H}_3$. For $i \in \{0,3\}$, let $H_i = \pi(\bar{H}_i)$. Clearly, $H_0$ and $H_3$ are $0$- and $3$-faces of $\B$, respectively. Then, $H_0 \cap H_3$ is a cycle of length $6$ of alternating colours $1$ and $2$ (each facet of $\B$ contains four such cycles; see Figure \ref{fig:vfX}). Moreover, only one link of $H_0 \cap H_3$ has non-trivial voltage and by Lemma \ref{lem:voltcon}, we see that $\pi^{-1}(H_0 \cap H_3)$ is connected and therefore a cycle of length $6n$ of alternating colours $1$ and $2$. Clearly, $\bar{H}_0 \cap \bar{H}_3 \subset \pi^{-1}(H_0 \cap H_3)$ and every connected component of $\bar{H}_0 \cap \bar{H}_3$ is a cycle of alternating colours $1$ and $2$. It follows that $\bar{H}_0 \cap \bar{H}_3$ has no other choice but to be the whole preimage $\pi^{-1}(H_0 \cap H_3)$. Therefore $\bigcap \mathcal{H} = \bar{H}_0 \cap \bar{H}_3 $ is a connected subgraph.

Now, suppose $\mathcal{H}$ contains an $i$-face, say $\bar{H}_{i_0}$, for some $i \in \{1,2\}$. Let $H_{i_j} = \pi(\bar{H}_{i_j})$ and let $\bar{\Delta} = \bigcap_{j = 0}^k \bar{H}_{i_j}$. Note that $\pi(\bar{\Delta})= \bigcap_{j = 0}^k H_{i_j}$. Since $H_{i_0}$ is an $i$-face with $i \in \{1,2\}$, $\pi^{-1}(H_{i_0})$ has $n$ connected components, one of which is $\bar{H}_{i_0}$. By Remark \ref{rem:iso} the restriction of $\pi$ to $\bar{H}_{i_0}$ is a graph isomorphism. In particular, since $\bar{\Delta} \subset \bar{H}_{i_0}$, the restriction $\pi \mid_{\bar{\Delta}}$ is a graph isomorphism between $\bar{\Delta}$ and $\bigcap_{j = 0}^k H_{i_j}$. However $\bigcap_{j = 0}^k H_{i_j}$ is connected, since $\B$ is a maniplex. Therefore $\bar{\Delta}$ is connected. We conclude that $\B^n$ has the CIP.

To see that the $0$- and $3$-faces of $\B^n$ are not bipartite, it suffices to show they contain an odd cycle. Consider a $3$-face $\bar{F}$ of $\B^n$ and note that $\bar{F}=\pi^{-1}(F)$ for some $3$-face $F$ of $\B$. Then, $F$ is isomorphic to the left hand side graph of Figure \ref{fig:vfX}. Let $u$ be a flag $2$-adjacent to a flag in $X$. The path $C$ starting at $u$ and tracing links of colours $012101021$ is a cycle of length $9$. Since all darts in $C$ have trivial voltage, $\pi^{-1}(C)$ consist of $n$ disjoint copies of $C$, all of which are contained in $\bar{F}$. Therefore $\bar{F}$ is not bipartite. The case when $\bar{F}$ is a $0$-face follows by duality.
\end{proof}


By combining the above lemma with Lemma \ref{lem:canCIP} we obtain the following.

\begin{proposition}
\label{prop:main}
For all $n \in \NN$, the canonical double cover $\overline{\B^n}$ of $\B^n$ is a faithful thin $4$-maniplex without the CIP.
\end{proposition}

\subsection{Regularity}

In this subsection we will show that the maniplexes $\overline{\B^n}$ of Proposition \ref{prop:main} are regular. By Remark \ref{rem:can}, it suffices to show that $\B^n$ is regular. That is, it is enough to show that for a fixed flag $(u,0)$ of $\B^n$ and for each $i \in [4]$, there exists an automorphism $\hat{\rho}_i \in \Aut(\B^n)$ mapping $(u,0)$ to its $i$-neighbour $(u,0)^i$. Three of the four automorphisms needed can be found by `lifting' a canonical set of generators of $\Aut(\B)$.

We say that an automorphism $\varphi$ of $\B$ {\em lifts} if there exists $\varphi' \in \Aut(\B^n)$ such that $\varphi' \circ \pi =  \pi \circ \varphi$ where $\pi \colon \B^n \to \B$ denotes the canonical projection. The automorphism $\varphi'$ is then called a {\em lift} of $\varphi$. Sufficient and necessary conditions for automorphisms of a graph to lift were first given in \cite{lift}. In particular, graph automorphisms that `commute' with an automorphism of the voltage group have a lift. The following lemma is (broadly speaking) a consequence of the results in \cite{lift} and is a particular case of \cite[Corollary 22]{GenVolt} (where we set $\omega(x)$ to be trivial everywhere).

\begin{lemma}
\label{lem:lift}
Let $(\Gamma,\zeta)$ be a voltage graph with voltage group $G$, and $\varphi \in \Aut(\Gamma)$ and $f \in \Aut(G)$. If $\zeta(x^\varphi) = f(\zeta(x))$ for all darts $x$ of $\Gamma$, then the permutation $\varphi'$ given by $(u,i) \mapsto (u^\varphi,f(i))$ is a lift of $\varphi$.
\end{lemma}

Recall that $\F_X$ is the set of flags of $\B$ that are incident to a dart in the set $X$, defined  in Section \ref{sec:covers}. That is $\F_X = \{u \in \B \mid uv \in X\}$. Let $u$ be a fixed flag in $\F_X$. For each $i \in [4]$, let $\rho_i$ be the automorphism of $\B$ mapping $u$ to $u^i$. We know such automorphisms exist because $\B$ is regular. By inspecting Figures \ref{fig:edgeX} and \ref{fig:vfX}, one can see that $\rho_0$ and $\rho_3$ preserve the voltages, while $\rho_2$ inverts them (this can also be readily checked by a computer program such as {\sf Gap}). However, $i \mapsto -i$ is a automorphism of the voltage group $\ZZ_n$. Therefore if $i \in \{0,2,3\}$, then Lemma \ref{lem:lift} tells us that $\rho_i$ has a lift $\rho_i' \in \Aut(\B^n)$. Note that $\rho_i'$ maps $(u,0)$ to some flag $(u^i,j)$, and since $\Aut(\B^n)$ acts transitively on the fibres, there exists $g \in \Aut(\B^n)$ mapping $(u^i,j)$ to $(u,0)^i$. Thus, the automorphism $\hat{\rho}_i := g \circ \rho_i'$ maps $(u_0,0)$ to its $i$-neighbour.

Now, it remains to show that there exists a $\hat{\rho}_1 \in \Aut(\B^n)$ such that $(u,0)^{\hat{\rho}_1} = (u,0)^1$. Lemma \ref{lem:lift} seems to be of little use here, so we will simply give $\hat{\rho}_1$ explicitly. Let $\F_X^{\rho_1} := \{ v^{\rho_1} \mid v \in \F_X\}$ and define $\hat{\rho}_1$ as the permutation given by  

\[
(v,i)^{\hat{\rho}_1}=
 \begin{cases}
(v^{\rho_1},-i) & \text{ if $v \in \F_X \cup \F_X^{\rho_1}$} \\ 
(v^{\rho_1},-i+1) & \text{ if $v \notin \F_X \cup \F_X^{\rho_1}$}

\end{cases}
\]

for all flags $(v,i)$ of $\B^n$. We must show that $\hat{\rho}_1$ is an automorphism of $\B^n$. Clearly $\hat{\rho}_1$ is bijective, as $\rho_1$ is an automorphism and both mapping $i \mapsto -i$ and $i \mapsto -i + 1$ are bijective. To prove that $\hat{\rho}_1$ preserves structure, we must show that if two flags of $\B^n$ are $i$-adjacent, then their images under $\hat{\rho}_1$ are $i$-adjacent or, equivalently, that $\hat{\rho}_1$ commutes with $\Mon(\B)$. This process is straightforward, but unfortunately rather tedious. Since the the image of a flag $(v,i)$ depends on where the flag $v$ is in $\B$, we have several cases to consider, depending on whether or not $v$ has an $i$-neighbour in $\F_X$ or in $\F_X^{\rho_1}$ for some $i \in [4]$. This amounts to nine different cases, but due to some of them being equivalent, we only need to analyse five of them. To better understand this, consider the following observations.



\begin{remark}
\label{rem:adjacencies}
If $v$ is a flag of $\B$ and let $j \in [4]$. Then the following hold:
\begin{enumerate}
\item for $j \in \{0,3\}$, $v \in \F_X$ if and only if $v^j \in \F_X$;
\item for $j \in \{1,2\}$, $v \in \F_X$ implies  $v^j \notin \F_X$;
\item $v^1 \in \F_X$ if and only if $v \in \F_X^{\rho_1}$;
\item if $j \neq 2$ then $(v,i)^j= (v^j,i)$;
\item if $j = 2$ then $(v,i)^j= (v^j,i+\alpha)$ where

\begin{itemize}
\item[(a)] $\alpha = 1$ if $v \in \F_X$,
\item[(b)] $\alpha = -1$ if $v^2 \in \F_X$,
\item[(c)] $\alpha = 0$ in all other cases.
\end{itemize}

\end{enumerate}
\end{remark}

By Remark \ref{rem:adjacencies}, the cases where $v^0 \in \F_X$, $v^3 \in \F_X$ and $v^1 \in \F_X^{\rho_1}$ are all equivalent to assuming $v \in \F_X$. Similarly the cases where $v^0 \in \F_X^{\rho_1}$, $v^3 \in \F_X^{\rho_1}$ and $v^1 \in \F_X$ are all equivalent to assuming $v \in \F_X^{\rho_1}$. Then we only need to consider the cases where $v \in \F_X$, $v \in \F_X^{\rho_1}$, $v^2 \in \F_X$, $v^2 \in \F_X^{\rho_1}$ and where none of the neighbours of $v$ are in $\F_X \cup \F_X^{\rho_1}$. We will do the full analysis for the first case, but omit the remaining four as they follow from a similar computation.

Let $v \in \F_X$, let $j \in [4]$ and notice that $v^j \in \F_X \cup \F_X^{\rho_1}$ only when $j \neq 2$, by Remark \ref{rem:adjacencies}. If $j \in \{0,1,3\}$, observe that 
\begin{align*}
 (v,i)^{j\hat{\rho}_1} = (v^j,i)^{\hat{\rho}_1} = (v^{j\rho_1},-i) = (v^{\rho_1j},-i) = (v^{\rho_1},-i)^j = (v,i)^{\hat{\rho}_1j}, 
\end{align*}
where the first and third equalities follow from Remark \ref{rem:adjacencies}, the second and fifth equalities from the definition of $\hat{\rho}_1$ and the third one from the fact that automorphisms and monodromy elements commute.
Similarly, if $j = 2$ then 
\begin{align*}
(v,i)^{j\hat{\rho}_1} = (v^j,i+1)^{\hat{\rho}_1} = (v^{j\rho_1},-(i+1)+1) = (v^{\rho_1j},-i) = (v^{\rho_1},-i)^j = (v,i)^{\hat{\rho}_1j}.
\end{align*}

This shows that $\bar{\rho}_1$ commutes with $\Mon(\B^n)$ when $v \in \F_X$. By a similar analysis for the remaining cases we obtain:


\begin{proposition}
\label{prop:rank4}
For all positive integers $n$, $\overline{\B^n}$ is regular.
\end{proposition}

\begin{figure}[h!]
\begin{center}
\includegraphics[width=0.70\textwidth]{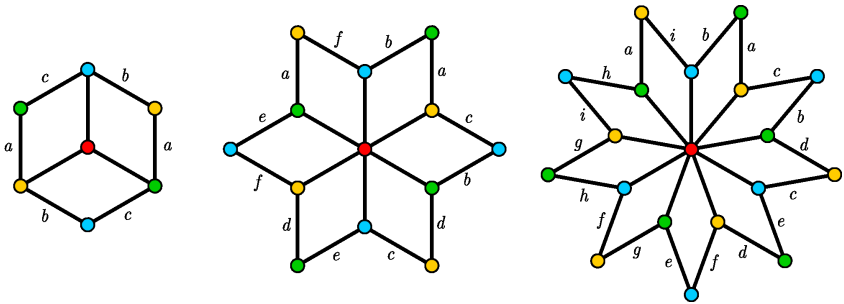}
\caption{From left to right, the skeleton of a facet of $\B$, $\B^2$ and $\B^3$, where vertices with the same colour are identified, as are edges with the same letter label. The facets and vertex figures of $\B^n$ are of Schl\"{a}fli type $\{4,3n\}$ and $\{3n,4\}$, respectively.}
\label{fig:facets}
\end{center}
\end{figure}



\section{Higher ranks}
\label{sec:higher}

So far we have proved that there are infinitely many regular $4$-maniplexes that are faithful and thin but non-polytopal. This result easily extends to every rank $n >4$ by constructing ``maniplex ravioli''. Given a rank $n$ maniplex $\M$, we can construct a rank $n+1$ maniplex by taking two copies of $\M$ and 'glueing' them in such a way that we obtain a maniplex having exactly two facets isomorphic to $\M$. The resulting maniplex will have $\Aut(\M) \times \ZZ_2$ as its automorphism group, and will be faithful, thin or polytopal if and only if $\M$ has the corresponding property as well. This construction is fairly well-known and can be seen as a trivial colour-coded extension (see \cite{colourcoded} for details).

\begin{definition}
Let $\M$ be an $n$-maniplex and let $\F(\M)$ denote its flag-set. The {\em raviolo} of $\M$, $R(\M)$, is the $(n+1)$-maniplex with flag-set $\F(\M) \times \ZZ_2$ and where the $i$-adjacencies are given by:
\[ (u,j)^i =
\begin{cases}
(u^i,j) & \text{ if $i \neq n$,}\\
(u,j+1) & \text{ if $i = n$.}\\
\end{cases}
\]
\end{definition}

For any two positive integers $n$ and $k$ we can apply the construction above $k$ times successively to $\overline{B^n}$ to obtain a regular $(4+k)$-maniplex that is faithful and thin but is non-polytopal. This completes the proof of Theorem \ref{the:main}.

\section{Concluding remarks}

The key for the construction presented in this paper to work is to carefully define voltage graphs (or voltage maniplexes, if you will) $(\M,\zeta)$ with the property that for any two faces of ranks $i-1$ and $i+1$, the image under $\pi^{-1}$ of one of those faces will have a single connected component, while the image of the other will have $|G|$ connected components (where $G$ is the voltage group). This property, which we may call {\em alternate connectedness}, is instrumental for properties such as faithfulness and thinness to be carried over from the base maniplex to its derived cover.

When $\zeta$ is the voltage assignment associated to the canonical double cover, it suffices that $\M$ is such that for any two faces of ranks $i-1$ and $i+1$, one of them is bipartite while the other is not, for $(\M,\zeta)$ to be alternately connected. We have chosen $\B$ to have this property and we have defined the voltages $\zeta_n$ in such a way that $(\B,\zeta_n)$ is alternately connected and the non-bipartiteness of the $0$- and $3$-faces is inherited by the cover $\B^n$.

Therefore one can obtain a thin, faithful non-polytopal $n$-maniplex for every polytopal $n$-maniplex whose faces are `alternately bipartite'. In rank $4$, the (flag-graphs of the) $11$-cell and $57$-cell are two well-known examples of polytopal self-dual maniplexes (with hemi-icosedral and hemi-dodecaheral facets, respectively) that have this property. Thus, their canonical double covers are thin and faithful but non-polytopal. Whether we can find voltage assignments for each of the $11$- and $57$-cell that are equivalent to the assignments $\zeta_n$ defined for $\B$ is not clear. The $57$-cell is too large to quickly check for regular $n$-fold covers in \cite{atlas} or {\sf Gap} (using the library of groups of small order). The $11$-cell has two regular (polytopal) double covers in the atlas \cite{atlas}, both of which have non-bipartite $0$- and $3$-faces, but no regular triple cover. However, regularity and polytopality are not necessary here; it suffices that the $11$-cell admits $n$-fold covers that are `alternately bipartite' for every $n>2$, even if it is not regular.




%

%

Regarding the construction in Section \ref{sec:higher}, the maniplexes of ranks $n > 4$ obtained might be seen as `degenerate', as they have only two facets. However, other methods for constructing extensions exist that yield maniplexes with more than two facets. Some non-trivial colour-coded extensions producing maniplexes with $4$ facets are known to preserve thinness and non-polytopality. Yet, the proof of this is rather lengthy and we leave it out of this paper for the sake of brevity.

\end{document}